\documentclass[1p]{article}

\usepackage{graphicx}
\usepackage{amsmath}
\usepackage{amssymb}
\usepackage{graphicx}
\usepackage{fancybox}
\usepackage{ascmac}
\usepackage{color}
\usepackage{amsthm}
\usepackage{amscd}
\usepackage{ulem}
\usepackage[margin=30mm]{geometry}

\newtheorem{thm}{Theorem}[section]
\newtheorem{dfn}[thm]{Definition}
\newtheorem{pro}[thm]{Proposition}
\newtheorem{lemm}[thm]{Lemma}

\newtheorem{prf}{Proof}

\begin{document}

\title{A new co-tame automorphism of the polynomial ring}
\author{Shoya Yasuda\footnote{Affiliation: Department of Mathematical Sciences, Tokyo Metropolitan University, 1-1 Minami-Osawa, Hachioji,
Tokyo 192-0397, Japan.}}
\date{}

\maketitle
\begin{abstract}
 In this paper, we discuss subgroups of the automorphism group of the polynomial ring in $n$ variables over a field of characteristic zero.
An automorphism $F$ is said to be {\it co-tame} if the subgroup generated by $F$ and affine automorphisms contains the tame subgroup. 
In 2017, Edo-Lewis gave a sufficient condition for co-tameness of automorphisms.
Let ${\it EL}_n$ be the set of all automorphisms satisfying Edo-Lewis's condition.
Then, for a certain topology on the automorphism group of the polynomial ring, any element of the closure of ${\it EL}_n$ is co-tame. 
Moreover, all the co-tame automorphisms previously known belong to the closure of ${\it EL}_n$. 
In this paper, we give the first example of co-tame automorphisms in $n$ variables which do not belong to the closure of ${\it EL}_n$.
\end{abstract}

\section{Introduction}

Let $k$ be a field of characteristic zero, $k[{\bf x}]:=k[x_1,\ldots,x_n]$ the polynomial ring in $n$ variables over $k$, and $\mathrm{Aut}_kk[{\bf x}]$ the automorphism group of the $k$-algebra $k[{\bf x}]$.
We write each $\phi \in \mathrm{Aut}_kk[{\bf x}]$ as $\phi=(\phi(x_1),\ldots,\phi(x_n))$.
We say that $\phi \in \mathrm{Aut}_kk[{\bf x}]$ is {\it affine} if $\phi = (x_1,\ldots,x_n)A+\mathbf{b}$ for some $A \in \mathit{GL}_n(k)$ and $\mathbf{b} \in k^n$, and {\it triangular} if 
$$\phi = (a_1 x_1 + f_1,\ldots, a_i x_i +f_i,\ldots, a_nx_n+f_n)$$
for some $a_i \in k^*$ and $f_i \in k[x_1,\ldots,x_{i-1}]$.
We denote by $\mathrm{Aff}_n(k)$ (resp.\ $\mathrm{BA}_n(k)$) the set of affine (resp.\ triangular) automorphisms of $k[{\bf x}]$.
Then, $\mathrm{Aff}_n(k)$ and $\mathrm{BA}_n(k)$ are the subgroups of $\mathrm{Aut}_kk[{\bf x}]$.
Note that $\tau_{\mathbf{a}}:=(x_1,\ldots,x_n)+\mathbf{a}$ is an element of $\mathrm{Aff}_n(k)$ for each $\mathbf{a} \in k^n$, and $\mathrm{Tr}_n(k):=\{ \tau_{\mathbf{a}} \mid \mathbf{a} \in k^n \}$ is a subgroup of $\mathrm{Aff}_n(k)$.

We say that $\phi \in \mathrm{Aut}_kk[{\bf x}]$ is {\it tame} if $\phi$ belongs to the {\it tame subgroup} $\mathrm{TA}_n(k):=\langle \mathrm{Aff}_n(k), \mathrm{BA}_n(k) \rangle$.
Jung \cite{Jung} and van der Kulk \cite{van1} showed that $\mathrm{Aut}_kk[x_1, x_2]=\mathrm{TA}_2(k)$.
In fact, $\mathrm{Aut}_kk[x_1, x_2]$ is the amalgamated free product of $\mathrm{Aff}_2(k)$ and $\mathrm{BA}_2(k)$ over $\mathrm{Aff}_2(k) \cap \mathrm{BA}_2(k)$ (cf.\ \cite[Part I, Theorem 3.3.]{Nagata}).
In 1972, Nagata \cite{Nagata} conjectured that $\psi \in \mathrm{Aut}_kk[x_1, x_2, x_3]$ defined by
\begin{align*}
\psi(x_1)=x_1+2(x_1x_3-x_2^2)x_2+(x_1x_3-x_2^2)^2x_3, \quad \psi(x_2)=x_2+(x_1x_3-x_2^2)x_3
\end{align*}
and $\psi(x_3)=x_3$ is not tame.
In 2004, this famous conjecture was solved in the affirmative by Shestakov-Umirbaev \cite{SU1}, \cite{SU2}.
When $n \geq 4$, it is not known whether $\mathrm{Aut}_kk[{\bf x}]=\mathrm{TA}_n(k)$.

In this paper, we study {\it co-tame} automorphisms defined as follows.
\begin{dfn}[Edo \cite{Edo1}] \rm
$\phi \in \mathrm{Aut}_kk[{\bf x}]$ is said to be {\it co-tame} if $\langle \phi,\mathrm{Aff}_n(k) \rangle \supset \mathrm{TA}_n(k)$.
\end{dfn}
No element of $\mathrm{Aut}_kk[x_1, x_2]$ is co-tame because of the amalgamated free product structure of $\mathrm{TA}_2(k)$ mentioned above (see \cite{EL}).
When $n \geq 3$, it is difficult to decide co-tameness of elements of $\mathrm{Aut}_kk[{\bf x}]$ in general.
The first example of co-tame automorphism was gave by Derksen.
He showed that the triangular automorphism $(x_1,\ldots,x_{n-1},x_n+x_1^2)$ is co-tame if $n \geq 3$ (cf.\ \cite[Theorem 5.2.1]{van}).
More generally, Bodnarchuck \cite{Bod1} showed that every non-affine element of $\mathrm{BA}_n(k) \circ \mathrm{Aff}_n(k) \circ \mathrm{BA}_n(k)$ is co-tame.
In 2013, Edo showed that a certain class of non-tame automorphisms, including Nagata's automorphism, are co-tame.
This result was recently generalized by Edo-Lewis \cite{EL2} as follows.
\begin{thm}[Edo-Lewis] \label{el}
Assume that $n\ge 3$. If $\phi \in \mathrm{Aut}_kk[x]\setminus \mathrm{Aff} _n(k)$ satisfies the following condition {\rm(}$\dag${\rm)}, then $\phi $ is co-tame{\rm:}\\
{\rm($\dag$) There exists ${\bf a}\in k^n\setminus \{ 0\}$ such that $\{ \phi \circ \tau _{\bf b}\circ \phi ^{-1}\mid {\bf b}\in k{\bf a}\} \subset \mathrm{Aff} _n(k)$}.
\end{thm}
We note that Edo-Lewis proved Theorem \ref{el} in the case where {\bf a} is a coordinate unit vector of $k^n$, but it easily implies the statement above.
Indeed, for any $\mathbf{a} \in k^n$, there exists $\alpha \in \mathrm{Aff}_n(k)$ such that 
$$\alpha \circ \{ \phi \circ \tau_{\mathbf{b}} \circ \phi^{-1} \mid \mathbf{b} \in k\mathbf{a} \} \circ \alpha^{-1}=\{ \psi \circ \tau_{\mathbf{c}} \circ \psi^{-1} \mid \mathbf{c} \in k(1,0,\ldots,0) \},$$
where $\psi:=\alpha \circ \phi \circ \alpha^{-1}$.
Moreover, $\psi$ is co-tame if and only if $\phi$ is co-tame.
In 2015, Edo-Lewis \cite{EL} gave the first example of automorphisms which are not co-tame.
They found such automorphisms in $\mathrm{TA}_3(k) \setminus \mathrm{Aff}_3(k)$.

The purpose of this paper is to construct a co-tame automorphism $\phi$ of $k[\mathbf{x}]$ with the following property:\\
($\ddag$) {\it $\langle \phi, \mathrm{Tr}_3(k) \rangle$ is the free product of $\langle \phi \rangle$ and $\mathrm{Tr}_n(k)$, and $\langle \phi, \mathrm{Tr}_n(k) \rangle \cap \mathrm{Aff}_n(k) = \mathrm{Tr}_n(k)$}.\\
Here, $\phi$ satisfies ($\ddag$) if and only if $\phi$ satisfies ${\phi}^{i_1} \circ \tau_{\mathbf{a}_1} \circ \dots \circ {\phi}^{i_{s-1}} \circ \tau_{\mathbf{a}_{s-1}} \circ {\phi}^{i_s} \notin \mathrm{Aff}_n(k)$ for any $i_1,\ldots,i_s \in \mathbb{Z} \setminus \{0\}$ and $\mathbf{a}_1,\ldots,\mathbf{a}_{s-1} \in k^n \setminus \{0\}$ with $s \geq 1$.

We claim that such co-tame automorphisms are essentially new.
To explain this, we introduce a topology which is due to S. Kuroda. In general, let $G$ be a group, and $H$ a subgroup of $G$. 
For each $S\subset G$, we define $\overline{S}^H$ to be the set of $g\in G$ such that $\langle g,H\rangle \cap S\ne \emptyset$.
Then, one can easily check that $2^G\ni S\mapsto \overline{S}^H\in 2^G$ is a closure operator. We remark that, if $A$ and $T$ are subsets of $G$ with $H\subset A$, then the set of $g\in G$ satisfying $\langle g,A\rangle \supset T$ is a closed subset of $G$ for this topology.

Now, we consider the topology on $\mathrm{Aut}_kk[\mathbf{x}]$ defined as above for $G=\mathrm{Aut}_kk[\mathbf{x}]$ and $H=\mathrm{Tr}_n(k)$.
Then, the set $\mathit{CT}_n$ of co-tame automorphisms of $k[\mathbf{x}]$ is a closed subset of $\mathrm{Aut} _kk[\mathbf{x}]$ by the remark.

Let $\mathit{EL}_n$ be the set of $\phi \in \mathrm{Aut}_kk[x]$ satisfying $(\dag $).
Then, we have
$$
\mathrm{Aff}_n(k) \circ \mathit{EL} _n \circ \mathrm{Aff}_n(k)\subset \mathit{EL}_n \subset \mathit{CT}_n.
$$
Since $\mathit{CT}_n$ is closed, we get $\overline{\mathit{EL}}_n^{\mathrm{Tr}_n(k)}\subset \mathit{CT}_n$.
It is notable that all the co-tame automorphisms previously known belong to $\overline{\mathit{EL}}_n^{\mathrm{Tr}_n(k)}$ (cf.\ \cite{Edo1}, \cite{EL2}). 
On the other hand, a co-tame automorphism satisfying ($\ddag$) does not belong to $\overline{\mathit{EL}}_n^{\mathrm{Tr}_n(k)}$.

\section{Main result}

A $k$-linear map $D:k[\mathbf{x}] \to k[\mathbf{x}]$ is called a {\it $k$-derivation} on $k[\mathbf{x}]$ if $D$ satisfies $D(fg)=fD(g)+D(f)g$ for all $f, g \in k[\mathbf{x}]$.
If $D$ is a $k$-derivation on $k[\mathbf{x}]$, then we can write
$$D=D(x_1)\frac{\partial}{\partial x_1} + \cdots + D(x_n) \frac{\partial}{\partial x_n}.$$
We say that a $k$-derivation $D$ on $k[\mathbf{x}]$ is {\it locally nilpotent} if there exists a positive integer $l$ such that $D^l(f)=0$ for every $f \in k[\mathbf{x}]$.
We denote by $\mathrm{Der}_kk[{\bf x}]$ (resp.\ $\mathrm{LND}_kk[{\bf x}]$) the set of $k$-derivations (resp.\ locally nilpotent $k$-derivations) on $k[\mathbf{x}]$.
If $D \in \mathrm{LND}_kk[{\bf x}]$ and $f \in \ker D$, we have $fD \in \mathrm{LND}_kk[{\bf x}]$ since $D(fg)=fD(g)$ for all $g \in k[\mathbf{x}]$.
We say that $D \in \mathrm{Der}_kk[{\bf x}]$ is {\it triangular} if $D(x_i) \in k[x_1,\ldots,x_{i-1}]$ for all $1 \leq i \leq n$.
If $D$ is triangular, then $D \in \mathrm{LND}_kk[{\bf x}]$. 
Given $D \in \mathrm{LND}_kk[{\bf x}]$, we define a {\it exponential automorphism} $\exp D \in \mathrm{Aut}_kk[{\bf x}]$ by
\begin{align*}
(\exp D)(f):=\sum_{i=0}^\infty \frac{D^i(f)}{i!} \quad \mathrm{for} \ \mathrm{all} \ f \in k[{\bf x}].
\end{align*}

We define $\sigma \in \mathrm{Aut}_kk[{\bf x}]$ by $\sigma(x_i):=x_{n-i+1}$ for $i=1,\ldots,n$, and set
\begin{align} 
D&:=\sum_{j=1}^{n-1} (n-j)x_{j+1} \frac{\partial}{\partial x_j}, \label{eq007} \\
D'&:=\sigma \circ D \circ \sigma = \sum_{j=2}^{n} (j-1)x_{j-1} \frac{\partial}{\partial x_j}. \label{eq008}
\end{align}
Since we can regard $D$ as a triangular derivation by changing the ordering of variables, $D$ and $D'$ are locally nilpotent.

We define a polynomial $f_{[n]}$ as follows.
When $n=2m-1$ with $m \geq 2$, we define
\begin{align} \label{eq009}
f_{[2m-1]}:=\frac{1}{2} \sum_{i=1}^{2m-1} (-1)^{i-1} \binom{2m-2}{i-1} x_{i}x_{2m-i}.
\end{align}
Then, we see that $\sigma(f_{[2m-1]})=f_{[2m-1]}$.
When $n=2m$ with $m \geq 2$, we define
\begin{align} \label{eq010}
f_{[2m]}:=D(f_{[2m-1]})^2-2D^2(f_{[2m-1]})f_{[2m-1]}.
\end{align}
We have the following theorem.
\begin{thm} \label{t1}
If $n \geq 3$, then $f_{[n]}$ belongs to $\ker D \cap \ker D'$.
\end{thm}
The following is our main theorem in this paper.
\begin{thm} \label{t}
In the notation above, the following assertions hold for each $l \geq 1${\rm :}\\
{\rm (i)} \ $\phi:=\exp f_{[n]}^lD \circ \exp D' \circ \exp (-f_{[n]}^lD)$ is co-tame.\\
{\rm (ii)} \ $G:=\langle \phi, \mathrm{Tr}_n(k) \rangle$ is the free product of $\langle \phi \rangle$ and $\mathrm{Tr}_n(k)$, and $G \cap \mathrm{Aff}_n(k)=\mathrm{Tr}_n(k)$.
Hence, $\phi$ satisfies {\rm(\ddag)}.
\end{thm}

\section{Proof of Theorem \ref{t1}}

\begin{flushleft}
(I) When $n=2m-1$ with $m \geq 2$
\end{flushleft}

Let $f:=f_{[2m-1]}$.
Here, we remark $\sigma(f)=f$.
Hence, if $f \in \ker D$, then we have $f \in \ker D'$ since $D':=\sigma \circ D \circ \sigma$.
Thus, it suffices to $f \in \ker D$.

By computation, we have
\begin{align*} 
D(2f)&=D \left( \sum_{i=1}^{2m-1} (-1)^{i-1} \binom{2m-2}{i-1} x_{i}x_{2m-i} \right) \\
&=\sum_{i=1}^{2m-1} (-1)^{i-1} \binom{2m-2}{i-1} D(x_{i})x_{2m-i} + \sum_{i=1}^{2m-1} (-1)^{i-1} \binom{2m-2}{i-1} x_{i}D(x_{2m-i}) \\
&=\sum_{i=1}^{2m-2} (-1)^{i-1} \binom{2m-2}{i-1} \cdot (2m-1-i) \cdot x_{i+1}x_{2m-i} \\ 
&\quad + \sum_{i=2}^{2m-1} (-1)^{i-1} \binom{2m-2}{i-1} \cdot (i-1) \cdot x_{i} x_{2m-i+1}.
\end{align*}
Since $\binom{2m-2}{i-1} \cdot (2m-1-i)=\binom{2m-2}{i} \cdot i$, the right-hand side of equality above is equal to
\begin{align*}
\sum_{i=1}^{2m-2} (-1)^{i-1} \binom{2m-2}{i} \cdot i \cdot x_{i+1}x_{2m-i} + \sum_{i=2}^{2m-1} (-1)^{i-1} \binom{2m-2}{i-1} \cdot (i-1) \cdot x_{i} x_{2m-i+1}=0.
\end{align*}
This proves $f \in \ker D$.

\begin{flushleft}
(II) When $n=2m$ with $m \geq 2$
\end{flushleft}

Let $f:=f_{[2m]}$ and $f':=f_{[2m-1]}$.
Then, we have
\begin{align*} 
D(f)&=D(D(f')^2-2D^2(f')f')\\
&=2D^2(f')D(f')-2D^3(f')f'-2D^2(f')D(f') \\
&=-2D^3(f')f'.
\end{align*}
Thus, we show $D^3(f')=0$.
By computation, we have
\begin{align*} 
&D(2f')\\
&\quad =D \left( \sum_{i=1}^{2m-1} (-1)^{i-1} \binom{2m-2}{i-1} x_{i}x_{2m-i} \right) \\
&\quad=\sum_{i=1}^{2m-1} (-1)^{i-1} \binom{2m-2}{i-1} D(x_{i})x_{2m-i} + \sum_{i=1}^{2m-1} (-1)^{i-1} \binom{2m-2}{i-1} x_{i}D(x_{2m-i}) \\
&\quad=\sum_{i=1}^{2m-1} (-1)^{i-1} \binom{2m-2}{i-1} \cdot (2m-i) \cdot x_{i+1}x_{2m-i} \\
&\quad \quad+ \sum_{i=1}^{2m-1}(-1)^{i-1} \binom{2m-2}{i-1} \cdot i \cdot x_{i}x_{2m-i+1}\\
&\quad=\sum_{i=2}^{2m} (-1)^{i-2} \binom{2m-2}{i-2} \cdot (2m-i+1) \cdot x_{i}x_{2m-i+1} \\
&\quad \quad+ \sum_{i=1}^{2m-1}(-1)^{i-1} \binom{2m-2}{i-1} \cdot i \cdot x_{i}x_{2m-i+1}\\
&\quad=2x_1x_{2m}+\sum_{i=2}^{2m-1} (-1)^{i-1} \left( \binom{2m-2}{i-1} \cdot i - \binom{2m-2}{i-2} \cdot (2m-i+1) \right)x_{i}x_{2m-i+1}.
\end{align*}
Since $\binom{2m-2}{i-2} \cdot (2m-i)=\binom{2m-2}{i-1} \cdot (i-1)$ for $i=2,\ldots,2m-1$, the right-hand side of equality above is equal to
\begin{align*} 
&2x_1x_{2m}+\sum_{i=2}^{2m-1} (-1)^{i-1} \left( \binom{2m-2}{i-1} \cdot i - \binom{2m-2}{i-1} \cdot (i-1) - \binom{2m-2}{i-2} \right)x_{i}x_{2m-i+1}\\
&\quad=2x_1x_{2m}+\sum_{i=2}^{2m-1} (-1)^{i-1} \left(\binom{2m-2}{i-1} - \binom{2m-2}{i-2} \right)x_{i}x_{2m-i+1}.
\end{align*}
Thus, we know that $\sigma(D(f'))=D(f')$.
Next, we calculate $D^2(2f')$.
For $i=2,\ldots 2m-1$, we have
\begin{align*} 
&D\left( \left(\binom{2m-2}{i-1} - \binom{2m-2}{i-2} \right)x_{i}x_{2m-i+1}\right)\\
&\quad=\left(\binom{2m-2}{i-1} - \binom{2m-2}{i-2}\right)(D(x_{i})x_{2m-i+1}+x_{i}D(x_{2m-i+1}))\\
&\quad=\left(\binom{2m-2}{i-1} - \binom{2m-2}{i-2}\right)((2m-i)x_{i+1}x_{2m-i+1}+(i-1)x_{i}x_{2m-i+2})\\
&\quad=\binom{2m-2}{i-1}(2m-2i+1)x_{i+1}x_{2m-i+1} +  \binom{2m-2}{i-2} (2m-2i+1) x_{i}x_{2m-i+2}.
\end{align*}
Thus, we have
\begin{align*} 
D^2(2f')&=2(2m-1)x_2x_{2m}+\sum_{i=2}^{2m-1}(-1)^{i-1} \binom{2m-2}{i-1}(2m-2i+1)x_{i+1}x_{2m-i+1}\\
&\quad + \sum_{i=2}^{2m-1}(-1)^{i-1}  \binom{2m-2}{i-2} (2m-2i+1) x_{i}x_{2m-i+2}\\
&=2(2m-1)x_2x_{2m}+\sum_{i=3}^{2m}(-1)^{i-2} \binom{2m-2}{i-2}(2m-2i+3)x_{i}x_{2m-i+2} \\
&\quad+ \sum_{i=2}^{2m-1}(-1)^{i-1}  \binom{2m-2}{i-2} (2m-2i+1) x_{i}x_{2m-i+2}\\
&=4x_2x_{2m}+2\sum_{i=3}^{2m-1}(-1)^{i-2} \binom{2m-2}{i-2}x_{i}x_{2m-i+2}\\
&=2\sum_{i=2}^{2m}(-1)^{i-2} \binom{2m-2}{i-2}x_{i}x_{2m-i+2}=4\sigma(f').
\end{align*}
Now, $D'(f')=\sigma \circ D \circ \sigma (f')=0$ holds by (I).
This implies that $D \circ \sigma (f')=0$.
Hence, we get $D^3(f')=0$.
Finally, since $$\sigma(f)=\sigma(D(f'))^2-2\sigma(D^2(f'))\sigma(f')=D(f')^2-2(2f')\left(\frac{1}{2}D^2(f')\right)=f,$$
we have $f \in \ker D \cap \ker D'$.

\section{Proof of Theorem \ref{t} (i)}
We fix an integer $l \geq 1$.
Let $D$, $D'$, $f:=f_{[n]}$ and $\phi$ be as in Theorem \ref{t}.
For $p \in \ker D$, $q \in \ker D'$ and $a, b \in k^*$, we set $d:=\deg(f) \in \{2, 4 \}$,
\begin{align*}
\epsilon_p &:= \exp pD, & \epsilon'_q&:= \exp qD',  \\
\quad \mu_a &:= (a^{n-1} x_1,\ldots, a^{n-2i+1}x_i ,\ldots, a^{-n+1}x_n), & \nu_b &:= (b x_1, \ldots, b x_{i} ,\ldots,b x_n).
\end{align*}
Then, $\phi$ can be written as $\epsilon_{f^l} \circ \epsilon'_1 \circ \epsilon_{-f^l}$.

\begin{lemm} \label{c}
The following equalities hold for any $a \in k^*$ and $b \in k$: \\
{\rm (i)} $\mu_a \circ \epsilon_{bf^l} = \epsilon_{a^{-2}bf^l} \circ \mu_{a}$. \\
{\rm (ii)} $\nu_a \circ \epsilon_{bf^l} = \epsilon_{a^{dl}bf^l} \circ \nu_a$. \\
{\rm (iii)} $\mu_{a} \circ \epsilon'_b = \epsilon'_{a^2b} \circ \mu_{a}$. \\
{\rm (iv)} $\nu_{a} \circ \epsilon'_{b} = \epsilon'_{b} \circ \nu_a$. 
\end{lemm}
\begin{prf} \rm
Since $\mu_a \circ D \circ \mu_a^{-1}=a^{-2}D$ and $\mu_a(f)=f$,
we have
$$\mu_a \circ bf^lD \circ \mu_a^{-1} = \mu_a(bf^l)\cdot( \mu_a \circ D \circ \mu_a^{-1})=a^{-2}bf^lD.$$
Since $\epsilon_{f^l}=\exp f^lD$, it follows that
$$\mu_a \circ \epsilon_{bf^l} \circ \mu_a^{-1} = \exp (\mu_a \circ bf^lD \circ \mu_a^{-1}) = \exp a^{-2}bf^lD = \epsilon_{a^{-2}bf^l}.$$
This proves (i).
We can prove (ii), (iii) and (iv) similarly.\\
\qed
\end{prf}

\begin{lemm} \label{d}
For each $a \in k^*$, we have $\epsilon_{1-a} \circ \epsilon'_{-1} \circ \epsilon_{1-a^{-1}} \circ \epsilon'_{a}=\mu_{a^{-1}}$.
\end{lemm}
\begin{prf} \rm
For $a \in k^*$, we define $2\times2$ matries $A(a)$, $B(a)$ and $C(a)$ by
$$
A(a):=\left(
    \begin{array}{cc}
      1 & 0  \\
      a & 1  
    \end{array}
  \right), \
B(a):=\left(
    \begin{array}{cc}
      1 & a  \\
      0 & 1 
    \end{array}
  \right) \ \mathrm{and} \
C(a):=\left(
    \begin{array}{cc}
      a & 0 \\
      0 & a^{-1} 
    \end{array}
  \right).
$$
Then, we have
\begin{align}\label{d1}
A(1-a)B(-1)A(1-a^{-1})B(a)=C(a^{-1}).
\end{align}

Let $V_{n-1}$ be a $k$-vector subspace of $k[x, y]$ generated by $\{ x^{n-i}y^{i-1} \ | \ i=1,\ldots,n \}$.
For $M(a) \in \{ A(a), B(a), C(a) \}$, we define a $k$-linear map $\sigma_{M(a)} : V_{n-1} \to V_{n-1}$ by $\sigma_{M(a)}=(x,y)M(a)$, and denote by $M'(a)$ the matrix representation of $\sigma_{M(a)}$ with reference to $\{ x^{n-i}y^{i-1} \ | \ i=1,\ldots,n \}$.
Then, we have 
\begin{align}\label{d2}
A'(1-a)B'(-1)A'(1-a^{-1})B'(a)=C'(a^{-1})
\end{align}
by (\ref{d1}).
We put
$$A'(a)=(\mathbf{a}_1(a),\ldots,\mathbf{a}_n(a)), \ B'(a)=(\mathbf{b}_1(a),\ldots,\mathbf{b}_n(a)) \ \mathrm{and} \ C'(a)=(\mathbf{c}_1(a),\ldots,\mathbf{c}_n(a)).$$
By computation, we get
$$\mathbf{a}_i(a)=\sum_{j=0}^{n-i} \binom{n-i}{j}a^{j}\mathbf{e}_{i+j}, \ \mathbf{b}_i(a)=\sum_{j=0}^{i-1} \binom{i-1}{j}a^{j}\mathbf{e}_{i-j} \ \mathrm{and} \ \mathbf{c}_i(a)=a^{n-2i+1}\mathbf{e}_{i}$$
for $i=1,\ldots,n$, where $\mathbf{e}_{1},\ldots,\mathbf{e}_{n}$ are the standard unit vectors of $k^n$.
Now, since
\begin{align*}
\epsilon_a(x_i)&=\sum_{j=0}^{\infty} \frac{a^j}{j!}{D}^j(x_i)=\sum_{j=0}^{n-i} \binom{n-i}{j} a^j x_{i+j},\\
\epsilon'_a(x_i)&=\sum_{j=0}^{\infty} \frac{a^j}{j!}(D')^j(x_i)=\sum_{j=0}^{i-1} \binom{i-1}{j} a^j x_{i-j}
\end{align*}
and $\mu_a(x_i)=a^{n-2i+1}x_{i}$ for $i=1,\ldots,n$, we have
$$\epsilon_a=(x_1,\ldots,x_n)A'(a), \ \epsilon'_a=(x_1,\ldots,x_n)B'(a) \ \mathrm{and} \ \mu_a=(x_1,\ldots,x_n)C'(a).$$
Thus by (\ref{d2}), we obtain $\epsilon_{1-a} \circ \epsilon'_{-1} \circ \epsilon_{1-a^{-1}} \circ \epsilon'_{a}=\mu_{a^{-1}}$.
\qed
\end{prf}

Now, we prove Theorem \ref{t} (i).
Fix $u \in k^*$ with $1-u^{2dl} \neq 0$, and define
\begin{align}
&\sigma_1:={\phi}^{-1} \circ \epsilon_{1-u^{-dl}} \circ \mu_{u^{2^{-1}dl}} \circ \nu_u \circ \phi, \nonumber \\
&\sigma_2:=\nu_{u^{-1}} \circ \mu_{u^{2^{-1}dl}} \circ \epsilon_{1-u^{dl}}. \nonumber
\end{align}

We show that $\sigma_1 \circ \sigma_2=\epsilon_{(1-u^{2dl})f^l}$.
By Lemmas \ref{c} and \ref{d}, we have
\begin{align}
\sigma_1&=\uwave{{\phi}^{-1}} \circ \epsilon_{1-u^{-dl}} \circ \mu_{u^{2^{-1}dl}} \circ \nu_u \circ \uwave{\phi} & & \nonumber \\
&=\epsilon_{f^l} \circ \epsilon'_{-1} \circ \epsilon_{-f^l} \circ \epsilon_{1-u^{-dl}} \circ \mu_{u^{2^{-1}dl}} \circ \uwave{\nu_u \circ \epsilon_{f^l}} \circ \epsilon'_1 \circ \epsilon_{-f^l} & & \rm{(Lemma \ \ref{c} (ii))} \nonumber \\
&=\epsilon_{f^l} \circ \epsilon'_{-1} \circ \epsilon_{-f^l} \circ \epsilon_{1-u^{-dl}} \circ \uwave{\mu_{u^{2^{-1}dl}} \circ \epsilon_{u^{dl}f^l}} \circ \nu_u \circ \epsilon'_1 \circ \epsilon_{-f^l} & & \rm{(Lemma \ \ref{c} (i))} \nonumber \\
&=\epsilon_{f^l} \circ \epsilon'_{-1} \circ \uwave{\epsilon_{-f^l}} \circ \epsilon_{1-u^{-dl}} \circ \uwave{\epsilon_{f^l}} \circ \mu_{u^{2^{-1}dl}} \circ \nu_u \circ \epsilon'_1 \circ \epsilon_{-f^l} \nonumber \\
&=\epsilon_{f^l} \circ \epsilon'_{-1} \circ \epsilon_{1-u^{-dl}} \circ \mu_{u^{2^{-1}dl}} \circ \uwave{\nu_u \circ \epsilon'_1} \circ \epsilon_{-f^l} & & \rm{(Lemma \ \ref{c} (iv))} \nonumber \\
&=\epsilon_{f^l} \circ \epsilon'_{-1} \circ \epsilon_{1-u^{-dl}} \circ \uwave{\mu_{u^{2^{-1}dl}} \circ  \epsilon'_1} \circ \nu_u \circ \epsilon_{-f^l} & &\rm{(Lemma \ \ref{c} (iii))} \nonumber \\
&=\epsilon_{f^l} \circ \uwave{\epsilon'_{-1} \circ \epsilon_{1-u^{-dl}} \circ \epsilon'_{u^{dl}}} \circ \mu_{u^{2^{-1}dl}}  \circ \nu_u \circ \epsilon_{-f^l} & & \rm{(Lemma \ \ref{d})} \nonumber \\
&=\epsilon_{f^l} \circ \epsilon_{u^{dl}-1} \circ \uwave{\mu_{u^{-dl}} \circ \mu_{u^{2^{-1}dl}}} \circ \nu_u \circ \epsilon_{-f^l}  \nonumber \\
&=\epsilon_{f^l} \circ \epsilon_{u^{dl}-1} \circ \mu_{u^{-2^{-1}dl}} \circ \uwave{\nu_u \circ \epsilon_{-f^l}} & & \rm{(Lemma \ \ref{c} (ii))} \nonumber \\
&=\epsilon_{f^l} \circ \epsilon_{u^{dl}-1} \circ \uwave{\mu_{u^{-2^{-1}dl}} \circ \epsilon_{-u^{dl}f^l}} \circ \nu_u & & \rm{(Lemma \ \ref{c} (i))} \nonumber \\
&=\uwave{\epsilon_{f^l}} \circ \epsilon_{u^{dl}-1} \circ \uwave{\epsilon_{-u^{2dl}f^l}} \circ \mu_{u^{-2^{-1}dl}} \circ \nu_u \nonumber \\
&=\epsilon_{(1-u^{2dl})f^l} \circ \epsilon_{u^{dl}-1} \circ \mu_{u^{-2^{-1}dl}} \circ \nu_u. \nonumber
\end{align}
Therefore, we get
\begin{align*}
\sigma_1 \circ \sigma_2 &= \epsilon_{(1-u^{2dl})f^l} \circ \uwave{\epsilon_{u^{dl}-1} \circ \mu_{u^{-2^{-1}dl}} \circ \nu_u \circ \nu_{u^{-1}} \circ \mu_{u^{2^{-1}dl}} \circ \epsilon_{1-u^{dl}}} = \epsilon_{(1-u^{2dl})f^l}.
\end{align*}
It is known that $\epsilon_{(1-u^{2dl})f^l}$ belongs to $\overline{\mathit{EL}}_n^{\mathrm{Tr}_n(k)}$  (cf.\ \cite[Theorem 12.]{EL2}).
Hence, $\epsilon_{(1-u^{2dl})f^l}$ is co-tame.
Since $\sigma_1 \circ \sigma_2$ belongs to $\langle \phi, \mathrm{Aff}_n(k) \rangle$, it follows that 
$$\mathrm{TA}_n(k) \subset \langle \epsilon_{(1-u^{2dl})f^l}, \mathrm{Aff}_n(k) \rangle \subset \langle \phi, \mathrm{Aff}_n(k) \rangle.$$
This proves that $\phi$ is co-tame.

\section{Proof of Theorem \ref{t} (ii)}

The goal of this section is to prove the following proposition, which implies Theorem~\ref{t} (ii).

\begin{pro}\label{pro001}
For any $i_1,\ldots,i_s \in \mathbb{Z} \setminus \{0\}$ and $\mathbf{a}_1,\ldots,\mathbf{a}_{s-1} \in k^n \setminus \{\mathbf{0}\}$ with $s \geq 1$, we have
\begin{align} \label{eq001}
\theta := {\phi}^{i_1} \circ \tau_{\mathbf{a}_1} \circ \dots \circ {\phi}^{i_{s-1}} \circ \tau_{\mathbf{a}_{s-1}} \circ {\phi}^{i_s} \notin \mathrm{Aff}_n(k).
\end{align}
\end{pro}

Let $k[t, {\bf x}]$ be the polynomial ring in $n+1$ variables over $k$, and $\pi : k[t, {\bf x}] \to k[{\bf x}]$ the substitution map defined by $t \mapsto f$.
We identify $k[t, {\bf x}]/(t-f) = k[t, {\bf x}]/\ker \pi$ with $k[{\bf x}]$ via the isomorphism induced from $\pi$.

Let $\widetilde{D}, \widetilde{D'} \in \mathrm{LND}_kk[t, {\bf x}]$ be the extensions of $D$ and $D'$ defined by $\widetilde{D}(t)=\widetilde{D'}(t)=0$, respectively.
For each $p \in \ker \widetilde{D}$ and $q \in \ker \widetilde{D'}$, we define $\widetilde{\epsilon}_p, \widetilde{\epsilon'_q}, \widetilde{\phi} \in \mathrm{Aut}_kk[t, {\bf x}]$ by
$$
\widetilde{\epsilon}_p := \exp p\widetilde{D}, \quad \widetilde{\epsilon'_q}:= \exp q\widetilde{D'} \quad {\rm and} \quad \widetilde{\phi}:=\widetilde{\epsilon}_{t^l} \circ \widetilde{\epsilon'_1} \circ \widetilde{\epsilon}_{-t^l}.
$$
Then, $\widetilde{\phi}$ fixes $t-f$, since $\widetilde{D}$ and $\widetilde{D'}$ kill $t-f$.
Hence, $\widetilde{\phi}$ induce an automorphism of $k[t, {\bf x}]/(t-f)$, which is equal to $\phi$ by construction.

For each $\mathbf{a}=(a_1, \ldots, a_n) \in k^n$, we define $\widetilde{\tau}_{\mathbf{a}} \in \mathrm{Aut}_kk[t, {\bf x}]$ by
$$\widetilde{\tau}_{\mathbf{a}}:=(t+\tau_{\mathbf{a}}(f)-f, x_1+a_1, \ldots , x_n+a_n)$$ 
Then, we have $\widetilde{\tau_{\mathbf{a}}}(t-f)=t-f$.
The automorphism of $k[{\bf x}]$ induced from $\widetilde{\tau}_{\mathbf{a}}$ is equal to $\tau_{\mathbf{a}}$.
Therefore, for $\theta$ in (\ref{eq001}), we have
\begin{align} \label{eq002} 
\theta(x_1) = \theta(\pi(x_1)) = \pi(\widetilde{\theta}(x_1)),
\end{align}
where $\widetilde{\theta} := {\widetilde{\phi}}^{i_1} \circ \widetilde{\tau}_{\mathbf{a}_1} \circ \dots \circ {\widetilde{\phi}}^{i_{s-1}} \circ \widetilde{\tau}_{\mathbf{a}_{s-1}} \circ {\widetilde{\phi}}^{i_s}$.

Next, for $p = {\displaystyle \sum_{i_0,\ldots,i_n \geq 0} u_{i_0,\ldots,i_n}t^{i_0}x_1^{i_1} \cdots x_n^{i_n}} \in k[t, \mathbf{x}] \setminus \{0\}$ with $u_{i_0,\ldots,i_n} \in k$, we define
$$\mathrm{supp}(p):=\{ (i_0,\ldots,i_n) \in \mathbb{N}^{n+1} \mid u_{i_0,\ldots,i_n} \neq 0 \}.$$
For $\mathbf{w}=(w_0,\ldots,w_n) \in \mathbb{N}^{n+1}$, we define
$$\deg_{\mathbf{w}}(p):=\max\{i_0w_0+\dots+i_nw_n \mid (i_0,\ldots,i_n) \in \mathrm{supp}(p)\}.$$
We write $\deg(p):=\deg_{(1,\ldots,1)}(p)$.
We denote by $\mathrm{lt}(p)$ the leading term of $p$ for the {\it lexicographic order}, i.e., the ordering defined by $t^{i_0}x_1^{i_1}\cdots x_n^{i_n} < t^{j_0}x_1^{j_1} \cdots x_n^{j_n}$ if $i_m < j_m$ for the first $m$ with $i_m \neq j_m$.

Now, we define $\mathbf{w}_1,\mathbf{w}_2,\mathbf{w}_3 \in \mathbb{N}^{n+1}$ as follows:
\begin{align*}
\mathbf{w}_1:=\sum_{i=1}^{n+1}\mathbf{e}_i, \quad 
\mathbf{w}_2:=(n-2)\mathbf{e}_1+\sum_{i=2}^{n+1}(2n-i)\mathbf{e}_i, \quad
\mathbf{w}_3:=\sum_{i=2}^{n+1}\mathbf{e}_i.
\end{align*}
For each $\alpha, \beta \geq 1$, we define:
\begin{align*}
&\mathcal{P}_{\alpha, \beta} := \{ p \in k[t, \mathbf{x}] \setminus \{0\} \mid \deg_{\mathbf{w}_1}(p) \leq \alpha+\beta, \deg_{\mathbf{w}_2}(p) \leq (n-2)\alpha+(n-1)\beta, \mathrm{lt}(p) \in k^*t^{\alpha}x_n^{\beta}  \}, \\
&\mathcal{Q}_{\alpha, \beta} := \{ p \in k[t, \mathbf{x}] \setminus \{0\} \mid \deg_{\mathbf{w}_3}(p) \leq \beta, \mathrm{lt}(p) \in k^*t^{\alpha}x_1^{\beta} \}, \\
&\mathcal{R}_{\alpha, \beta} := \{ p \in k[t, \mathbf{x}] \setminus \{0\} \mid \deg_{\mathbf{w}_3}(p) \leq \beta, \mathrm{lt}(p) \in k^*t^{\alpha}x_n^{\beta} \}. \\
\end{align*}

We show that $\widetilde{\phi}^d(x_i)$ belongs to $\mathcal{P}_{(2n-i-1)l, 1}$ for any $d \in \mathbb{Z} \setminus \{0\}$ and $i=1,\ldots,n$.
We have
\begin{align*}
\widetilde{\epsilon}_{-t^l}(x_i)=\sum_{j=0}^{\infty} \frac{(-t^l)^j}{j!}\widetilde{D}^j(x_i)=\sum_{j=0}^{n-i} \binom{n-i}{j} (-t^l)^j x_{i+j}.
\end{align*}
We can easily check that $\widetilde{\epsilon}_{-t^l}(x_i) \in \mathcal{R}_{(n-i)l, 1}$.

\begin{lemm}\label{lem001}
For each $\alpha, \beta \geq 1$ and $d \in \mathbb{Z} \setminus \{0\}$ , we have $\widetilde{\epsilon'_d}(\mathcal{R}_{\alpha, \beta}) \subset \mathcal{Q}_{\alpha, \beta}$.
\end{lemm}
\begin{prf} \rm
For $i=1,\ldots,n$, we set $X_i:=\widetilde{\epsilon'_d}(x_i)$ and $T:=\widetilde{\epsilon'_d}(t)$.
Then, we have
$$X_i=\sum_{j=0}^{\infty} \frac{d^j}{j!}\widetilde{D'}^j(x_i)=\sum_{j=0}^{i-1} \binom{i-1}{j} d^j x_{i-j}$$
and $T=t$.
Thus, we have $\deg_{\mathbf{w}_3}(X_i)=1$ and $\deg_{\mathbf{w}_3}(T)=0$.

Take any $r \in \mathcal{R}_{\alpha, \beta}$.
For $\mathbf{i}=(i_0,\ldots,i_n) \in \mathrm{supp}(r)$, we set
$$d_{\mathbf{i}, \mathbf{w}_3} := \deg_{\mathbf{w}_3}(\widetilde{\epsilon'_d}(t^{i_0}x_1^{i_1} \cdots x_n^{i_n}))=\deg_{\mathbf{w}_3}(T^{i_0}X_1^{i_1}\cdots X_n^{i_n}).$$
Since $i_1+\cdots+i_n \leq \beta$, we get 
\begin{align*}
d_{\mathbf{i}, \mathbf{w}_3}=i_1+\cdots+i_n \leq \beta.
\end{align*}
This prove that $\deg_{\mathbf{w}_3}(r) \leq \beta$.
Similarly, we have 
$$\mathrm{lt}(\widetilde{\epsilon'_d}(t^{i_0}x_1^{i_1} \cdots x_n^{i_n}))=\mathrm{lt}(T^{i_0}X_1^{i_1} \cdots X_n^{i_n}) \in k^*m_\mathbf{i},$$
where $m_\mathbf{i} := t^{i_0}x_1^{i_1+\cdots+i_n}$.
If $\mathbf{i}=(\alpha,0,\ldots,0,\beta)$, then $m_\mathbf{i}=t^{\alpha}x_1^{\beta}$.
If $\mathbf{i}=(\alpha,0,\ldots,0,i_n)$ and $i_n < \beta$, then $m_\mathbf{i}=t^{\alpha}x_1^{i_n}$ is less than $t^{\alpha}x_1^{\beta}$.
If $\mathbf{i}=(i_0,\ldots,i_n) \in \mathrm{supp}(r)$ satisfies $i_0 < \alpha$, then $m_\mathbf{i}$ is less than $t^{\alpha}x_1^{\beta}$.
Thus, we get $\mathrm{lt}(\widetilde{\epsilon'_d}(r)) \in k^*t^{\alpha}x_1^{\beta}$.
Therefore, $\widetilde{\epsilon'_d}(r)$ belongs to $\mathcal{Q}_{\alpha, \beta}$.
\qed
\end{prf}

\begin{lemm} \label{lem002}
For each $\alpha, \beta \geq 1$, we have $\widetilde{\epsilon}_{t^l}(\mathcal{Q}_{\alpha, \beta}) \subset \mathcal{P}_{\alpha', \beta'}$, where $\alpha' := \alpha + (n-1)l\beta$ and $\beta' := \beta$.
\end{lemm}
\begin{prf} \rm
For $i=1,\ldots,n$, we set $X_i:=\widetilde{\epsilon}_{t^l}(x_i)$ and $T:=\widetilde{\epsilon}_{t^l}(t)$.
Then, we have
$$X_i=\sum_{j=0}^{\infty} \frac{(t^l)^j}{j!}\widetilde{D}^j(x_i)=\sum_{j=0}^{n-i} \binom{n-i}{j} (t^l)^j x_{i+j}$$
and $T=t$.
Thus, we have $\deg_{\mathbf{w}_1}(X_i)=(n-i)l+1$, $\deg_{\mathbf{w}_2}(X_i)=(n-2)(n-i)l+(n-1)$, $\deg_{\mathbf{w}_1}(T)=1$ and $\deg_{\mathbf{w}_2}(T)=n-2$.

Take any $q \in \mathcal{Q}_{\alpha, \beta}$.
For $\mathbf{i}=(i_0,\ldots,i_n) \in \mathrm{supp}(q)$ and $\mathbf{w} \in \{ \mathbf{w}_1, \mathbf{w}_2 \}$, we set
$$d_{\mathbf{i}, \mathbf{w}} := \deg_{\mathbf{w}}(\widetilde{\epsilon'_d}(t^{i_0}x_1^{i_1} \cdots x_n^{i_n}))=\deg_{\mathbf{w}}(T^{i_0}X_1^{i_1} \cdots X_n^{i_n}).$$
Since $i_0 \leq \alpha$ and $i_1+\cdots+i_n \leq \beta$, we get 
\begin{align*}
d_{\mathbf{i}, \mathbf{w}_1}&=i_0+\sum_{j=1}^n ((n-j)l+1) \cdot i_j =i_0 + \sum_{j=1}^n (n-j)l \cdot i_j + \sum_{j=1}^n i_j \\
&\leq\alpha + (n-1)l\beta +\beta = \alpha'+\beta' ,\\
d_{\mathbf{i}, \mathbf{w}_2}&=(n-2)i_0+\sum_{j=1}^n ((n-2)(n-j)l+(n-1)) \cdot i_j \\
&=(n-2)i_0 + \sum_{j=1}^n (n-2)(n-j)l \cdot i_j +\sum_{j=1}^n (n-1) \cdot i_j\\
&\leq (n-2)\alpha + (n-2)(n-1)l\beta +(n-1)\beta = (n-2)\alpha'+(n-1)\beta'.
\end{align*}
This prove that $\deg_{\mathbf{w}_1}(q) \leq \alpha'+\beta'$ and $\deg_{\mathbf{w}_2}(q) \leq (n-2)\alpha'+(n-1)\beta'$.
Similarly, we have
$$\mathrm{lt}(\widetilde{\epsilon}_{t^l}(t^{i_0}x_1^{i_1} \cdots x_n^{i_n}))=\mathrm{lt}(T^{i_0}X_1^{i_1} \cdots X_n^{i_n}) \in k^*m_\mathbf{i},$$
where $m_\mathbf{i} := t^{i_0+((n-1)i_1+(n-2)i_2+\cdots+i_{n-1})l}x_n^{i_1+\cdots+i_n}$.
If $\mathbf{i}=(\alpha,\beta,0,\ldots,0)$, then $m_\mathbf{i}=t^{\alpha'}x_n^{\beta'}$.
If $\mathbf{i}=(\alpha,i_1,\ldots,i_n) \in \mathrm{supp}(q)$ and $i_1 < \beta$, then 
\begin{align*}
i_0+((n-1)i_1+(n-2)i_2+\ldots+i_{n-1})l &\leq \alpha+(i_1+(n-2)(i_1+\cdots+i_n))l \\
&\leq \alpha+(i_1+(n-2)\beta)l \\
&< \alpha+(n-1)l\beta.
\end{align*}
Hence, $m_\mathbf{i}$ is less than $t^{\alpha'}x_n^{\beta'}$.
If $\mathbf{i}=(i_0,\ldots,i_n) \in \mathrm{supp}(q)$ satisfies $i_0 < \alpha$, then $m_\mathbf{i}$ is less than $t^{\alpha'}x_n^{\beta'}$ because $i_0+((n-1)i_1+(n-2)i_2+\ldots+i_{n-1})l \leq i_0+(n-1)l\beta < \alpha+(n-1)l\beta$.
Thus, we get $\mathrm{lt}(\widetilde{\epsilon}_{t^l}(q)) \in k^*t^{\alpha'}x_n^{\beta'}$.
Therefore, $\widetilde{\epsilon}_{t^l}(q)$ belongs to $\mathcal{P}_{\alpha', \beta'}$.
\qed
\end{prf}
By Lemmas \ref{lem001} and \ref{lem002}, we have $\widetilde{\phi}^d(x_i) \in \mathcal{P}_{(2n-i-1)l, 1}$ for each $d \in \mathbb{Z} \setminus \{0\}$ and $i=1,\ldots,n$.
We remark that, if $p$ is an element of $\mathcal{P}_{\alpha, \beta}$, then one of the following holds for each $\mathbf{i}=(i_0,\ldots,i_n) \in \mathrm{supp}(p-\mathrm{lt}(p))$:

(a) $\mathbf{i}=(\alpha,0,\ldots,0,i_n)$ and $0 \leq i_n < \beta$

(b) $0 \leq i_0 < \alpha$, ${\displaystyle \sum_{j=0}^{n} i_j} \leq \alpha+\beta$ and $(n-2)i_0+{\displaystyle \sum_{j=1}^{n}(2n-j-1) \cdot i_j} \leq (n-2)\alpha+(n-1)\beta$

\begin{lemm} \label{lem003}
For each $p \in \mathcal{P}_{\alpha, \beta}$, we have $\deg(\pi(p)) = \deg(f)\alpha+\beta \geq 3$.
\end{lemm}
\begin{prf} \rm
Note that $d(\mathbf{i}):=\deg(\pi(t^{i_0}x_1^{i_1} \cdots x_n^{i_n}))=\deg(f^{i_0}x_1^{i_1} \cdots x_n^{i_n})=\deg(f)i_0+i_1+\cdots+i_n$ for each $\mathbf{i}=(i_0,\ldots,i_n) \in \mathrm{supp}(p)$.
If $\mathbf{i}$ is as in (a) or (b), then $d(\mathbf{i})$ is less than $d((\alpha,0,\ldots,0,\beta))=\deg(f)\alpha+\beta$.
Hence, $\deg(\pi(p))$ is equal to $\deg(\pi(\mathrm{lt}(p)))=\deg(\pi(t^{\alpha}x_n^{\beta}))=\deg(f)\alpha+\beta$.
\qed
\end{prf}

In the following, we show that $\widetilde{\theta}(x_1)$ belongs to $\mathcal{P}_{\alpha, \beta}$ for some $\alpha, \beta \geq 1$.
This implies that $\theta \notin \mathrm{Aff}_n(k)$ by (\ref{eq002}) and Lemma \ref{lem003}.

\begin{flushleft}
(iii) When $n$ is an odd number
\end{flushleft}

Let $n$ be an odd number.
Take any $\mathbf{a}=(a_1, \ldots, a_n) \in k^n \setminus \{\mathbf{0}\}$, and set $\widetilde{\psi}:=\widetilde{\phi}^d \circ \widetilde{\tau}_{\mathbf{a}}$,
\begin{align*}
X_i := \widetilde{\psi}(x_i) = \widetilde{\phi}^d(x_i)+a_i \ \mathrm{for} \ i=1,\ldots,n \ \mathrm{and} \ T:= \widetilde{\psi}(t).
\end{align*}
By (\ref{eq009}), we have
$$T \in k^*t+k^*\sum_{j=1}^{n} a_{n+1-j}X_j + k.$$
Then, we have the following, where $\mu:=\max\{i \mid a_i \neq 0\}+n-2$:
\begin{align*}
\deg_{\mathbf{w}_1}(X_i)&=(2n-i-1)l+1, & \deg_{\mathbf{w}_1}(T)&=\mu l+1, \\
\deg_{\mathbf{w}_2}(X_i)&=(n-2)(2n-i-1)l+(n-1), & \deg_{\mathbf{w}_2}(T)&=(n-2)\mu l+(n-1),\\
\mathrm{lt}(X_i)&\in k^*t^{(2n-i-1)l} x_n, & \mathrm{lt}(T) &\in k^* t^{\mu l} x_n.
\end{align*}

Now, observe that $\widetilde{\phi}^d(x_1)$ belongs to $\mathcal{P}_{(2n-2)l,1}$.
Hence, we obtain from the following proposition that $\widetilde{\theta}(x_1) \in \mathcal{P}_{\alpha, \beta}$ for some $\alpha, \beta \geq 1$ as claimed.

\begin{pro}
For each $\alpha, \beta \geq 1$, we have $\widetilde{\psi}(\mathcal{P}_{\alpha, \beta}) \subset \mathcal{P}_{\alpha', \beta'}$, where $\alpha' := l(\mu \alpha + (n-1)\beta)$ and $\beta' := \alpha+ \beta$.
\end{pro}
\begin{prf} \rm
Take any $p \in \mathcal{P}_{\alpha, \beta}$.
For $\mathbf{i}=(i_0,\ldots,i_{n}) \in \mathrm{supp}(p)$ and $\mathbf{w} \in \{\mathbf{w}_1, \mathbf{w}_2\}$, we set
$$d_{\mathbf{i}, \mathbf{w}} := \deg_{\mathbf{w}}(\widetilde{\psi}(t^{i_0}x_1^{i_1} \cdots x_n^{i_n}))=\deg_{\mathbf{w}}(T^{i_0}X_1^{i_1} \cdots X_n^{i_n}).$$
Note that ${\displaystyle \sum_{j=0}^n i_j} \leq \alpha+\beta$ and $(n-2)i_0+{\displaystyle \sum_{j=1}^n (2n-j-1)\cdot i_j} \leq (n-2)\alpha+(n-1)\beta$.
Since $\mathrm{lt}(p) \in k^*t^{\alpha}x_n^{\beta}$ and $\mu > n-2$, we also have $i_0 \leq \alpha$ and $\mu i_0+{\displaystyle \sum_{j=1}^n (2n-j-1)\cdot i_j} \leq \mu \alpha+(n-1)\beta$.
Hence, we get 
\begin{align*}
d_{\mathbf{i}, \mathbf{w}_1}&=(\mu l + 1)i_0 + \sum_{j=1}^{n} ((2n-j-1)l+1) \cdot i_j\\
&= l \left(\mu i_0+\sum_{j=1}^{n} (2n-j-1) \cdot i_j \right)+ \sum_{j=0}^n i_j\\
&\leq l(\mu \alpha + (n-1)\beta) + \alpha + \beta =\alpha' + \beta', \\
d_{\mathbf{i}, \mathbf{w}_2}&= ((n-2)\mu l + (n-1)) \cdot i_0 + \sum_{j=1}^n ((n-2)(2n-j-1)l+(n-1)) \cdot i_j \\
&=(n-2)l \left( \mu  i_0 + \sum_{j=1}^{n} (2n-j-1) \cdot i_j \right) + (n-1)\sum_{j=0}^n i_j \\
&\leq (n-2)l(\mu \alpha + (n-1)\beta) + (n-1)(\alpha + \beta) = (n-2)\alpha' + (n-1)\beta'.
\end{align*}
This proves that $\deg_{\mathbf{w}_1}(p) \leq \alpha' + \beta'$ and $\deg_{\mathbf{w}_2}(p) \leq (n-2)\alpha' + (n-1)\beta'$.
Similarly, we have $\mathrm{lt}(\widetilde{\psi}(t^{i_0}x_1^{i_1} \cdots x_n^{i_n}))=\mathrm{lt}(T^{i_0}X_1^{i_1} \cdots X_n^{i_n}) \in k^*m_\mathbf{i}$, where
$$m_\mathbf{i} := t^{l(\mu i_0+(2n-2)i_1+(2n-3)i_2+\cdots+(n-1)i_n)}x_n^{i_0+\cdots+i_n}.$$
If $\mathbf{i}=(\alpha,0,\ldots,0,\beta)$, then $m_\mathbf{i}=t^{\alpha'}x_n^{\beta'}$.
If $\mathbf{i}$ is as in (a), then $m_\mathbf{i}=t^{l(\mu \alpha+(n-1)i_n)}x_n^{\alpha+i_n}$ is less than $t^{\alpha'}x_n^{\beta'}$ because $i_n < \beta$.
If $i$ is as in (b), then $l(\mu i_0+(2n-2)i_1+(2n-3)i_2+\cdots+(n-1)i_n) < l(\mu\alpha+(n-1)\beta)=\alpha'$.
Hence, $m_\mathbf{i}$ is less than $t^{\alpha'}x_n^{\beta'}$.
Thus, we get $\mathrm{lt}(\widetilde{\psi}(p)) \in k^*t^{\alpha'}x_n^{\beta'}$.
Therefore, $\widetilde{\psi}(p)$ belongs to $\mathcal{P}_{\alpha', \beta'}$.
\qed
\end{prf}

\begin{flushleft}
(iv) When $n$ is an even number
\end{flushleft}

Let $n:=2m$ with $m \geq 2$.
Take any $\mathbf{a}=(a_1, \ldots, a_{2m}) \in k^{2m} \setminus \{\mathbf{0}\}$, and set $\widetilde{\psi}:=\widetilde{\phi}^d \circ \widetilde{\tau}_{\mathbf{a}}$,
\begin{align*}
X_i := \widetilde{\psi}(x_i) = \widetilde{\phi}^d(x_i)+a_i \ \mathrm{for} \ i=1,\ldots,2m \ \mathrm{and} \ T:= \widetilde{\psi}(t).
\end{align*}
Similar to (iii), we determine $\deg_{\mathbf{w}_1}(T)$, $\deg_{\mathbf{w}_2}(T)$ and $\mathrm{lt}(T)$.
We set $f:=f_{[2m]}$ and $f':=f_{[2m-1]}$.
By a simple calculating, we have
\begin{align*}
\tau_{\mathbf{a}}(f') &= f'+\sum_{i=2}^{2m}(-1)^{i} \binom{2m-2}{i-2} \cdot a_{2m-i+1}x_{i-1} + p(a_1,\ldots,a_{2m-1}),\\
\tau_{\mathbf{a}}(D(f')) &= D(f')+a_{2m}x_1+ \sum_{i=2}^{2m-1}(-1)^{i-1} \left(\binom{2m-2}{i-1}-\binom{2m-2}{i-2} \right)a_{2m-i+1}x_i \\
&\quad+ a_1x_{2m}+q(a_1,\ldots,a_{2m}),\\
\tau_{\mathbf{a}}(D^2(f')) &= D^2(f') + 2\sum_{i=1}^{2m-1}(-1)^{i-1} \binom{2m-2}{i-1} \cdot a_{2m-i+1}x_{i+1}+r(a_2,\ldots,a_{2m}),
\end{align*}
where
\begin{align*}
p(a_1,\ldots,a_{2m-1}) &\in k^*\sum_{i=1}^{m}a_i \cdot a_{2m-i},\\
q(a_1,\ldots,a_{2m}) &\in k^*\sum_{i=1}^{m}a_i \cdot a_{2m-i+1}, \\
r(a_2,\ldots,a_{2m}) &\in k^*\sum_{i=1}^{m}a_{i+1} \cdot a_{2m-i+1}.
\end{align*}
Here, we note that $\widetilde{\tau}_{\mathbf{a}}(t)=t+\tau_{a}(f)-f$ and $f=D(f')^2-2D^2(f')f'$.
Hence, the total degree of $\widetilde{\tau}_{\mathbf{a}}(t)$ is three.
We write $\widetilde{\tau}_{\mathbf{a}}(t)=U_0+U_1+U_2+U_3$, where $U_i$ is a homogeneous polynomial of degree i.
Then, $U_3$ is
\begin{align}
&2\sum_{i=1}^{2m-1}(-1)^{i-1}\binom{2m-2}{i-1} \cdot a_{2m-i+1} (D(f')x_{i}-2f'x_{i+1}) \nonumber\\
&+2\sum_{i=2}^{2m}(-1)^{i-1} \binom{2m-2}{i-2} \cdot a_{2m-i+1}(D^2(f')x_{i-1}-D(f')x_{i}) \label{eq15}.
\end{align}
In order to determine the degree of $\widetilde{\phi}^d(U_3)$, we calculate $\widetilde{\phi}^d(D(f')x_{i}-2f'x_{i+1})$ and $\widetilde{\phi}^d(D^2(f')x_{i-1}-D(f')x_{i})$.
For $p \in \{ t^l, -t^l \}$, we have
\begin{align*}
\widetilde{\epsilon}_{p}(f'x_{i+1})&=\sum_{j=0}^{2m-i-1}p^j \cdot \binom{2m-i-1}{j} \cdot f'x_{i+j+1} + \sum_{j=1}^{2m-i}p^j \cdot \binom{2m-i-1}{j-1} \cdot D(f')x_{i+j}\\
&\quad+\frac{1}{2}\sum_{j=2}^{2m-i+1}p^j \cdot \binom{2m-i-1}{j-2} \cdot D^2(f')x_{i+j-1},\\
\widetilde{\epsilon}_{p}(D(f')x_{i})&=\sum_{j=0}^{2m-i}p^j \cdot \binom{2m-i}{j} \cdot D(f')x_{i+j} + \sum_{j=1}^{2m-i+1}p^j \cdot \binom{2m-i}{j-1} \cdot D^2(f')x_{i+j-1},\\
\widetilde{\epsilon}_{p}(D^2(f')x_{i-1})&=\sum_{j=0}^{2m-i+1}p^j \cdot \binom{2m-i+1}{j} \cdot D^2(f')x_{i+j-1}.
\end{align*}
Since $\binom{a}{b}-\binom{a-1}{b-1}=\binom{a-1}{b}$, we have
\begin{align}
\widetilde{\epsilon}_{p}(D(f')x_{i}-2f'x_{i+1})&=\sum_{j=0}^{2m-i-1}p^j \cdot \binom{2m-i-1}{j} \cdot (D(f')x_{i+j}-2f'x_{i+j-1}) \nonumber \\
&\quad+\sum_{j=1}^{2m-i}p^j \cdot \binom{2m-i-1}{j-1} \cdot (D^2(f')x_{i+j-1}-D(f')x_{i+j}), \label{eq11}\\
\widetilde{\epsilon}_{p}(D^2(f')x_{i-1}-D(f')x_{i})&=\sum_{j=0}^{2m-i}p^j \cdot \binom{2m-i}{j} \cdot (D^2(f')x_{i+j-1}-D(f')x_{i+j}).\label{eq12}
\end{align}
Here, we note that the polynomial $D^2(f')x_{2m}$ dose not appear in the right-hand sides of (\ref{eq11}) and (\ref{eq12}).
From the definition of $D'$, we see that $D'(f')=0$, $D'(D(f'))=2f'$ and $D'(D^2(f'))=2D(f')$.
Thus, for each $d \in \mathbb{Z} \setminus \{0\}$, we have
\begin{align*}
\widetilde{\epsilon'_d}(f'x_{i+1})&=\sum_{j=0}^{i}d^j \cdot \binom{i}{j} \cdot f'x_{i-j+1},\\
\widetilde{\epsilon'_d}(D(f')x_{i})&=2\sum_{j=1}^{i}d^j \cdot \binom{i-1}{j-1} \cdot f'x_{i-j+1} + \sum_{j=0}^{i-1}d^j \cdot \binom{i-1}{j} \cdot D(f')x_{i-j},\\
\widetilde{\epsilon'_d}(D^2(f')x_{i-1})&=2\sum_{j=2}^{i}d^j \cdot \binom{i-2}{j-2} \cdot f'x_{i-j+1} + 2\sum_{j=1}^{i-1}d^j \cdot \binom{i-2}{j-1} \cdot D(f')x_{i-j}\\
& \quad+\sum_{j=0}^{i-2}d^j \cdot \binom{i-2}{j} \cdot D^2(f')x_{i-j-1}.
\end{align*}
Since $\binom{a}{b}-\binom{a-1}{b-1}=\binom{a-1}{b}$, we have
\begin{align}
\widetilde{\epsilon'_d}(D(f')x_{i}-2f'x_{i+1})&=\sum_{j=0}^{i-1}d^j \cdot \binom{i-1}{j} \cdot (D(f')x_{i-j}-2f'x_{i-j+1}),\label{eq13}\\
\widetilde{\epsilon'_d}(D^2(f')x_{i-1}-D(f')x_{i})&=\sum_{j=1}^{i-1}d^j \cdot \binom{i-2}{j-1} \cdot (D(f')x_{i-j}-2f'x_{i-j+1})\nonumber \\
&\quad+\sum_{j=0}^{i-2}d^j \cdot \binom{i-2}{j} \cdot (D^2(f')x_{i-j-1}-D(f')x_{i-j})\label{eq14}.
\end{align}
Here, we note that the polynomial $f'x_{1}$ dose not appear in the right-hand sides of (\ref{eq13}) and (\ref{eq14}).
By (\ref{eq15}) through (\ref{eq14}), for $\mu:=\max\{i \mid a_i \neq 0\}$, we get
\begin{align}
\mathrm{lt}(\widetilde{\phi}^d(U_3)) &\leq t^{(2m+\mu-2)l}(D^2(f')x_{2m-1}-D(f')x_{2m}),\nonumber\\
\deg_{\mathbf{w}_1}(\widetilde{\phi}^d(U_3)) &\leq (2m+\mu-2)l+3, \label{eq100} \\
\deg_{\mathbf{w}_2}(\widetilde{\phi}^d(U_3)) &\leq (2m-2)(2m+\mu-2)l+8m-4.\nonumber
\end{align}

In the following, we determine the degree of $\widetilde{\phi}^d(\widetilde{\tau}_{\mathbf{a}}(t)-U_3)$.
In considering the degree of $\widetilde{\phi}^d(\widetilde{\tau}_{\mathbf{a}}(t)-U_3)$, the main term of $\widetilde{\tau}_{\mathbf{a}}(t)-U_3$ is
\begin{align*}
a_{2m}^2x_1^2+\sum_{i=2}^{2m-1}a_{2m-i+1}^2(c_{1,i} \cdot x_{i}^2+c_{2,i} \cdot x_{i-1}x_{i+1})+a_1^2x_{2m}^2,
\end{align*}
where $c_{1,i}:=\left(\binom{2m-2}{i-1}-\binom{2m-2}{i-2}\right)^2$ and $c_{2,i}:=4\left(\binom{2m-2}{i-1} \cdot \binom{2m-2}{i-2}\right)$.
We note that $c_{1,i}, c_{2,i} > 0$ for any $i=2,\ldots,2m-1$.
Now, we have $\widetilde{\phi}^d(x_i) \in \mathcal{P}_{(4m-i-1)l, 1}$ for any $i=1,\ldots,2m$.
Moreover, since the signs of $\mathrm{lt}(\widetilde{\phi}^d(x_{i-1}))$ and $\mathrm{lt}(\widetilde{\phi}^d(x_{i+1}))$ are the  same, $\mathrm{lt}(\widetilde{\phi}^d(c_{1,i} \cdot x_{i}^2+c_{2,i} \cdot x_{i-1}x_{i+1}))$ belongs to $k^*\mathrm{lt}(\widetilde{\phi}^d(x_{i}^2))$.
Thus, for $\mu:=\max\{i \mid a_i \neq 0\}$, we get
\begin{align}
\mathrm{lt}(\widetilde{\phi}^d(\widetilde{\tau}_{\mathbf{a}}(t)-U_3)) &\in k^* t^{2(2m+\mu-2)l}x_{2m}^2,\nonumber \\
\deg_{\mathbf{w}_1}(\widetilde{\phi}^d(\widetilde{\tau}_{\mathbf{a}}(t)-U_3)) &= 2(2m+\mu-2)l+2,\label{eq101}\\
\deg_{\mathbf{w}_2}(\widetilde{\phi}^d(\widetilde{\tau}_{\mathbf{a}}(t)-U_3)) &= 2(2m-2)(2m+\mu-2)l+2(2m-1).\nonumber
\end{align}
Set $\mu:=2\max\{i \mid a_i \neq 0\}+4m-4$.
Then, we have the following by (\ref{eq100}) and (\ref{eq101}):
\begin{align*}
\deg_{\mathbf{w}_1}(X_i)&=(4m-i-1)l+1, & \deg_{\mathbf{w}_1}(T)&=\mu l+2, \\
\deg_{\mathbf{w}_2}(X_i)&=(2m-2)(4m-i-1)l+(2m-1), & \deg_{\mathbf{w}_2}(T)&=(2m-2)\mu l+2(2m-1),\\
\mathrm{lt}(X_i)&\in k^*t^{(4m-i-1)l} x_{2m}, & \mathrm{lt}(T) &\in k^* t^{\mu l} x_{2m}^2.
\end{align*}

Now, observe that $\widetilde{\phi}^d(x_1)$ belongs to $\mathcal{P}_{(4m-2)l,1}$.
Hence, the following proposition implies that $\widetilde{\theta}(x_1) \in \mathcal{P}_{\alpha, \beta}$ for some $\alpha, \beta \geq 1$ as claimed.

\begin{pro}
For each $\alpha, \beta \geq 1$, we have $\widetilde{\psi}(\mathcal{P}_{\alpha, \beta}) \subset \mathcal{P}_{\alpha', \beta'}$, where $\alpha' := l(\mu \alpha + (2m-1)\beta)$ and $\beta' := 2\alpha+ \beta$.
\end{pro}
\begin{prf} \rm
Take any $p \in \mathcal{P}_{\alpha, \beta}$.
For $\mathbf{i}=(i_0,\ldots,i_{2m}) \in \mathrm{supp}(p)$ and $\mathbf{w} \in \{\mathbf{w}_1, \mathbf{w}_2\}$, we set
$$d_{\mathbf{i}, \mathbf{w}} := \deg_{\mathbf{w}}(\widetilde{\psi}(t^{i_0}x_1^{i_1} \cdots x_{2m}^{i_{2m}}))=\deg_{\mathbf{w}}(T^{i_0}X_1^{i_1} \cdots X_{2m}^{i_{2m}})$$
and $\mu:=2\max\{ i \mid a_i \neq 0\}+4m-4$.
Note that $2i_0+{\displaystyle \sum_{j=1}^{2m} i_j} \leq 2\alpha+\beta$ and $$\mu i_0+{\displaystyle \sum_{j=1}^{2m}(4m-i-1)i_j}  \leq \mu \alpha+(2m-1)\beta.$$
Since $\mathrm{lt}(p) \in k^*t^{\alpha}x_n^{\beta}$, we also have $i_0 \leq \alpha$.
Hence, we get 
\begin{align*}
d_{\mathbf{i}, \mathbf{w}_1}&=(\mu l+2)i_0+\sum_{j=1}^{2m}((4m-i-1)l+1)i_j \\
&=l \left( \mu i_0 + \sum_{j=1}^{2m}(4m-i-1)i_j  \right) + 2i_0 + \sum_{j=1}^{2m}i_j \\
&\leq l(\mu \alpha + (2m-1)\beta) + 2\alpha + \beta =\alpha' + \beta', \\
d_{\mathbf{i}, \mathbf{w}_2}&=((2m-2)\mu l+2(2m-1))i_0 + \sum_{j=1}^{2m} ((2m-2)(4m-i-1)l+(2m-1))i_j \\
&=(2m-2)l \left( \mu i_0 + \sum_{j=1}^{2m}(4m-i-1)i_j  \right) + (2m-1) \left( 2i_0 + \sum_{i=1}^{2m}i_j \right) \\
&\leq (2m-2)l(\mu \alpha + (2m-1)\beta) + (2m-1)(2\alpha + \beta) = (2m-2)\alpha' + (2m-1)\beta'.
\end{align*}
This proves that $\deg_{\mathbf{w}_1}(p) \leq \alpha' + \beta'$ and $\deg_{\mathbf{w}_2}(p) \leq (2m-2)\alpha' + (2m-1)\beta'$.
Similarly, we have $\mathrm{lt}(\widetilde{\psi}(t^{i_0}x_1^{i_1} \cdots x_{2m}^{i_{2m}}))=\mathrm{lt}(T^{i_0}X_1^{i_1} \cdots X_{2m}^{i_{2m}}) \in k^*m_\mathbf{i}$, where 
$$m_\mathbf{i} := t^{l(\mu i_0+(4m-2)i_1+(4m-3)i_2+\cdots+(2m-1)i_{2m})}x_{2m}^{2i_0+i_1+\cdots+i_{2m}}.$$
If $\mathbf{i}=(\alpha,0,\ldots,0,\beta)$, then $m_\mathbf{i}=t^{\alpha'}x_{2m}^{\beta'}$.
If $\mathbf{i}$ is as in (a), then $m_\mathbf{i}=t^{l(\mu \alpha+(2m-1)i_{2m})}x_{2m}^{\alpha+i_{2m}}$ is less than $t^{\alpha'}x_n^{\beta'}$, since $i_{2m} < \beta$.
If $\mathbf{i}$ is as in (b), then 
\begin{align*}
&l \left( \mu i_0 + \sum_{j=1}^{2m}(4m-i-1)i_j  \right) \\
&\quad=l \left((\mu-(2m-2)) i_0+ (2m-2) i_0 + \sum_{j=1}^{2m}(4m-i-1)i_j \right)\\
&\quad\leq l((\mu-(2m-2))i_0 + (2m-2)\alpha+(2m-1)\beta)\\
&\quad<  l(\mu\alpha+(2m-1)\beta)=\alpha'.
\end{align*}
Hence, $m_\mathbf{i}$ is less than $t^{\alpha'}x_{2m}^{\beta'}$.
Thus, we get $\mathrm{lt}(\widetilde{\psi}(p)) \in k^*t^{\alpha'}x_{2m}^{\beta'}$.
Therefore, $\widetilde{\psi}(p)$ belongs to $\mathcal{P}_{\alpha', \beta'}$.
\qed
\end{prf}

\section*{Acknowledgment}

The author thanks Professor Sigeru Kuroda and Professor Ryuji Tanimoto for useful discussions. 
This work was supported by JSPS KAKENHI Grant Number JP19J20334.


\begin{thebibliography}{20}
\bibitem{Bod1} Yu.\ V. Bodnarchuk, Generating properties of triangular and bitriangular birational automorphisms of an affine space, Dopov. Nats. Akad. Nauk Ukr. Mat. Prirodozn. Tekh. Nauki 2002, no.~11, 7--12.
\bibitem{Edo1} E. Edo, Coordinates of $R[x,y]$: constructions and classifications, Comm. Algebra {\bf 41} (2013), no.~12, 4694--4710.
\bibitem{EL} E. Edo\ and\ D. Lewis, The affine automorphism group of $\Bbb A^3$ is not a maximal subgroup of the tame automorphism group, Michigan Math. J. {\bf 64} (2015), no.~3, 555--568.
\bibitem{EL2} E. Edo\ and\ D. Lewis, Co-tame polynomial automorphisms, Internat. J. Algebra Comput. {\bf 29} (2019), no.~5, 803--825.
\bibitem{van} A. van den Essen, Polynomial Automorphisms and the Jacobi Conjecture, Progr. Math., {\bf 190}, Birkh\"{a}user, Basel, 2000.
\bibitem{Jung} H. W. E. Jung, \"Uber ganze birationale Transformationen der Ebene, J. Reine Angew. Math. {\bf 184} (1942), 161--174.
\bibitem{van1} W. van der Kulk, On polynomial rings in two variables, Nieuw Arch. Wiskunde (3) {\bf 1} (1953), 33--41.
\bibitem{Nagata} M. Nagata, On Automorphism Group of $k[x,y]$, Lectures in Mathematics, Department of Mathematics, Kyoto University, Vol. 5, Kinokuniya Book-Store Co. Ltd., Tokyo, 1972.
\bibitem{SU1} I. P. Shestakov\ and\ U. U. Umirbaev, Poisson brackets and two-generated subalgebras of rings of polynomials, J. Amer. Math. Soc. {\bf 17} (2004), no.~1, 181--196.
\bibitem{SU2} I. P. Shestakov\ and\ U. U. Umirbaev, The tame and the wild automorphisms of polynomial rings in three variables, J. Amer. Math. Soc. {\bf 17} (2004), no.~1, 197--227.
\end{thebibliography}
\end{document}